\documentclass[11pt]{amsart}
\usepackage{amssymb}
\usepackage{amsmath}
\usepackage[active]{srcltx}
\usepackage{t1enc}
\usepackage[latin2]{inputenc}
\usepackage{verbatim}
\usepackage{amsmath,amsfonts,amssymb,amsthm}
\usepackage[mathcal]{eucal}
\usepackage{enumerate}
\usepackage[centertags]{amsmath}
\usepackage{graphics}

\newtheorem{theorem}{Theorem}

\begin{document}
\author{George Tephnadze}
\title[Fejér means]{A note on the norm convergence by Vilenkin-Fejér means}
\address{G. Tephnadze, Department of Mathematics, Faculty of Exact and
Natural Sciences, Iv. Javaxikhishvili, Tbilisi State University,
Chavchavadze str. 1, Tbilisi 0128, Georgia}
\thanks{The research was supported by Shota Rustaveli National Science
Foundation grant no.52/54 (Bounded operators on the martingale Hardy spaces).%
}
\email{giorgitephnadze@gmail.com}
\date{}
\maketitle

\begin{abstract}
The main aim of this paper is to find necessary and sufficient conditions
for the convergence of Fejér means in terms of the modulus of continuity on
the Hardy spaces $H_{p},$ when $0<p\leq 1/2.$
\end{abstract}



Let $\mathbb{N}_{+}$ denote the set of the positive integers, $\mathbb{N}:=%
\mathbb{N}_{+}\cup \{0\}.$

Let $m:=(m_{0,}m_{1},...)$ denote a sequence of the positive integers not
less than 2.

Denote by 
\begin{equation*}
Z_{m_{k}}:=\{0,1,...,m_{k}-1\}
\end{equation*}%
the additive group of integers modulo $m_{k}.$

Define the group $G_{m}$ as the complete direct product of the group $%
Z_{m_{j}}$ with the product of the discrete topologies of $Z_{m_{j}}$`s.

The direct product $\mu $ of the measures 
\begin{equation*}
\mu _{k}\left( \{j\}\right) :=1/m_{k}\text{ \qquad }(j\in Z_{m_{k}})
\end{equation*}
is the Haar measure on $G_{m_{\text{ }}}$with $\mu \left( G_{m}\right) =1.$

If the sequence $m:=(m_{0,}m_{1},...)$ is bounded than $G_{m}$ is called a
bounded Vilenkin group, else it is called an unbounded one. \textbf{In this
paper we discuss bounded Vilenkin groups only.}

The elements of $G_{m}$ represented by sequences 
\begin{equation*}
x:=(x_{0},x_{1},...,x_{j},...)\qquad \left( \text{ }x_{k}\in
Z_{m_{k}}\right) .
\end{equation*}

It is easy to give a base for the neighborhood of $G_{m}$ 
\begin{equation*}
I_{0}\left( x\right) :=G_{m},
\end{equation*}%
\begin{equation*}
I_{n}(x):=\{y\in G_{m}\mid y_{0}=x_{0},...,y_{n-1}=x_{n-1}\}\text{ }(x\in
G_{m},\text{ }n\in \mathbb{N}).
\end{equation*}%
Denote $I_{n}:=I_{n}\left( 0\right) ,$ for $n\in \mathbb{N}$ and $\overline{%
I_{n}}:=G_{m}$ $\backslash $ $I_{n}$. Set 
\begin{equation*}
e_{n}:=\left( 0,...,0,1,0,...\right) \in G_{m},
\end{equation*}%
the $n$-th coordinate of which is 1 and the rest are zeros $\left( n\in 
\mathbb{N}\right) .$

If we define the so-called generalized number system based on $m$ as $\
M_{0}:=1,$ $M_{k+1}:=m_{k}M_{k\text{ }}\ (k\in \mathbb{N})$ \ then every $%
n\in \mathbb{N}$ can be uniquely expressed as $n=\sum_{k=0}^{\infty
}n_{j}M_{j},$ where $n_{j}\in Z_{m_{j}}$ $~(j\in \mathbb{N})$ and only a
finite number of $n_{j}$`s differ from zero. Let $\left\vert n\right\vert
:=\max $ $\{j\in \mathbb{N},$ $n_{j}\neq 0\}.$

Denote by $L_{1}\left( G_{m}\right) $ the usual (one dimensional) Lebesgue
space.

Define the complex valued function $r_{k}\left( x\right) :G_{m}\rightarrow 
\mathbb{C}
,$ the generalized Rademacher functions as 
\begin{equation*}
r_{k}\left( x\right) :=\exp \left( 2\pi \iota x_{k}/m_{k}\right) \text{
\qquad }\left( \iota ^{2}=-1,\text{ }x\in G_{m},\text{ }k\in \mathbb{N}%
\right)
\end{equation*}

Now define the Vilenkin system $\psi :=(\psi _{n}:n\in \mathbb{N})$ on $%
G_{m} $ as: 
\begin{equation*}
\psi _{n}(x):=\overset{\infty }{\underset{k=0}{\Pi }}r_{k}^{n_{k}}\left(
x\right) \text{ \qquad }\left( n\in \mathbb{N}\right)
\end{equation*}

Specifically, we call this system the Walsh-Paley one if $m\equiv 2$. The
Vilenkin system is orthonormal and complete in $L_{2}\left( G_{m}\right) \,$%
\cite{AVD,Vi}.

If $f\in L_{1}\left( G_{m}\right) $ we can establish the the Fourier
coefficients, the partial sums of the Fourier series, the Fejér means, the
Dirichlet and Fejér kernels in the usual manner: 
\begin{eqnarray*}
\widehat{f}\left( k\right) &:&=\int_{G}f\overline{\psi _{k}}d\mu \,,\text{ \
\ \ \ \ \ }\left( \text{ }n\in \mathbb{N}\right) \\
S_{n}f &:&=\sum_{k=0}^{n-1}\widehat{f}\left( k\right) \psi _{k}\text{ },%
\text{ \ \ }\left( n\in \mathbb{N}_{+},S_{0}f:=0\right) , \\
\sigma _{n}f &:&=\frac{1}{n}\sum_{k=1}^{n}S_{k}f,\text{ \ }%
D_{n}:=\sum_{k=0}^{n-1}\psi _{k\text{ }},\,\text{\ \ }\,K_{n}:=\frac{1}{n}%
\overset{n}{\underset{k=1}{\sum }}D_{k},\text{ \ \thinspace }\left( \text{ }%
n\in \mathbb{N}_{+}\text{ }\right) .
\end{eqnarray*}

Recall that 
\begin{equation}
\quad \hspace*{0in}D_{M_{n}}\left( x\right) =\left\{ 
\begin{array}{l}
M_{n},\text{\thinspace\ \thinspace \thinspace \thinspace if\thinspace
\thinspace }x\in I_{n} \\ 
0,\text{\thinspace \thinspace \thinspace\ \ \ \thinspace \thinspace if
\thinspace \thinspace }x\notin I_{n}%
\end{array}%
\right.  \label{3}
\end{equation}%
\vspace{0pt}

Let $q_{A}=M_{2A}+M_{2A-2}+...+M_{2}+M_{0},$ $\ 2<A\in \mathbb{N}_{+}$. Then
(see \cite{BGG,GoAMH}) 
\begin{equation}
q_{A-1}\left\vert K_{q_{A-1}}(x)\right\vert \geq \frac{M_{2k}M_{2s}}{4},%
\text{ for }x\in I_{2s+1}\left( s_{2k}e_{2k}+s_{2l}e_{2l}\right) ,
\label{3a}
\end{equation}%
where $1\leq s_{2k}\leq m_{2k}-1$ and $1\leq s_{2s}\leq m_{2s}-1,$ $%
k=0,1,...,A-3,$ $s=k+2,k+3,...,A-1.$

The norm (or quasinorm) of the space $L_{p}(G_{m})$ is defined by \qquad
\qquad \thinspace\ 
\begin{equation*}
\left\Vert f\right\Vert _{p}:=\left( \int_{G_{m}}\left\vert f\right\vert
^{p}d\mu \right) ^{1/p}\qquad \left( 0<p<\infty \right) .
\end{equation*}%
The space $L_{p,\infty }\left( G_{m}\right) $ consists of all measurable
functions $f$ for which

\begin{equation*}
\left\Vert f\right\Vert _{L_{p},\infty }:=\underset{\lambda >0}{\sup }%
\lambda ^{p}\mu \left( f>\lambda \right) <+\infty
\end{equation*}

The $\sigma $-algebra generated by the intervals $\left\{ I_{n}\left(
x\right) :x\in G_{m}\right\} $ will be denoted by $\digamma _{n}$ $\left(
n\in \mathbb{N}\right) .$ Denote by $f=\left( f^{\left( n\right) },n\in 
\mathbb{N}\right) $ a martingale with respect to $\digamma _{n}$ $\left(
n\in \mathbb{N}\right) .$ (for details see e.g. \cite{We1}). The maximal
function of a martingale $f$ is defend by \qquad 
\begin{equation*}
f^{\ast }=\sup_{n\in \mathbb{N}}\left\vert f^{\left( n\right) }\right\vert .
\end{equation*}

For $0<p<\infty $ the Hardy martingale spaces $H_{p}$ $\left( G_{m}\right) $
consist of all martingales, for which 
\begin{equation*}
\left\Vert f\right\Vert _{H_{p}}:=\left\Vert f^{\ast }\right\Vert
_{p}<\infty .
\end{equation*}

If $f=\left( f^{\left( n\right) },\text{ }n\in \mathbb{N}\right) $ is
martingale then the Vilenkin-Fourier coefficients must be defined in a
slightly different manner: $\qquad \qquad $ 
\begin{equation*}
\widehat{f}\left( i\right) :=\lim_{k\rightarrow \infty
}\int_{G_{m}}f^{\left( k\right) }\left( x\right) \overline{\Psi }_{i}\left(
x\right) d\mu \left( x\right) .
\end{equation*}

The concept of modulus of continuity in $H_{p}$ $\left( 0<p\leq 1\right) $
is defined by 
\begin{equation*}
\omega \left( 1/M_{n},f\right) _{H_{p}}:=\left\Vert f-S_{M_{n}}f\right\Vert
_{H_{p}}.
\end{equation*}%
\qquad \qquad \qquad \qquad

A bounded measurable function $a$ is p-atom, if there exist a dyadic
interval $I,$ such that%
\begin{equation*}
\int_{I}ad\mu =0,\text{ \ \ }\left\Vert a\right\Vert _{\infty }\leq \mu
\left( I\right) ^{-1/p},\text{ \ \ supp}\left( a\right) \subset I.\qquad
\end{equation*}

The dyadic Hardy martingale spaces $H_{p}$ $\left( G\right) $ for $0<p\leq 1$
have an atomic characterization. Namely the following theorem is true (see 
\cite{We3}):

\textbf{Theorem W. }A martingale $f=\left( f_{n},\text{ }n\in \mathbb{N}%
\right) $ is in $H_{p}\left( 0<p\leq 1\right) $ if and only if there exists
a sequence $\left( a_{k},\text{ }k\in \mathbb{N}\right) $ of p-atoms and a
sequence $\left( \mu _{k},\text{ }k\in \mathbb{N}\right) $ of a real numbers
such that for every $n\in \mathbb{N}$

\begin{equation}
\qquad \sum_{k=0}^{\infty }\mu _{k}S_{M_{n}}a_{k}=f_{n}  \label{2A}
\end{equation}%
and

\begin{equation*}
\qquad \sum_{k=0}^{\infty }\left\vert \mu _{k}\right\vert ^{p}<\infty ,
\end{equation*}%
Moreover, $\left\Vert f\right\Vert _{H_{p}}\backsim \inf \left(
\sum_{k=0}^{\infty }\left\vert \mu _{k}\right\vert ^{p}\right) ^{1/p},$
where the infimum is taken over all decomposition of\textit{\ }$f$\textit{\ }%
of the form (\ref{2A}).

The weak type-$\left( 1,1\right) $ inequality for maximal operator of Fejér
means $\sigma ^{\ast }$ can be found in Schipp \cite{Sc} for Walsh series
and in Pál, Simon \cite{PS} for bounded Vilenkin series. Fujji \cite{Fu} and
Simon \cite{Si2} verified that $\sigma ^{\ast }$ is bounded from $H_{1}$ to $%
L_{1}$. Weisz \cite{We2} generalized this result and proved the boundedness
of $\sigma ^{\ast }$ from the martingale space $H_{p}$ to the space $L_{p}$
for $p>1/2$. Simon \cite{Si1} gave a counterexample, which shows that
boundedness does not hold for $0<p<1/2.$ The counterexample for $p=1/2$ due
to Goginava \cite{Go}, (see also \cite{BGG2} and \cite{gog3}). Weisz \cite%
{we4} proved that $\sigma ^{\ast }$ is bounded from the Hardy space $H_{1/2}$
to the space $L_{1/2,\infty }$.

The author in \cite{tep2} and \cite{tep3} (for Walsh system see \cite%
{GoSzeged}) proved that the maximal operator $\widetilde{\sigma }_{p}^{\ast
} $ with respect Vilenkin system defined by 
\begin{equation*}
\widetilde{\sigma }_{p}^{\ast }:=\sup_{n\in \mathbb{N}}\frac{\left\vert
\sigma _{n}\right\vert }{\left( n+1\right) ^{1/p-2}\log ^{2\left[ 1/2+p%
\right] }\left( n+1\right) },
\end{equation*}%
where $0<p\leq 1/2$ and $\left[ 1/2+p\right] $ denotes integer part of $%
1/2+p,$ is bounded from the Hardy space $H_{p}$ to the space $L_{p}.$
Moreover, we showed that the order of deviant behaviour of the $n$-th Fejér
means was given exactly. As a corollary we get 
\begin{equation}
\left\Vert \sigma _{n}f\right\Vert _{p}\leq c_{p}\left( n+1\right)
^{1/p-2}\log ^{2\left[ 1/2+p\right] }\left( n+1\right) \left\Vert
f\right\Vert _{H_{p}}.  \label{1}
\end{equation}

For Walsh-Kaczmarz system analogical theorems are proved in \cite{GNCz} and 
\cite{tep4}.

Móricz and Siddiqi \cite{Mor} investigates the approximation properties of
some special Nörlund means of Walsh-Fourier series of $L_{p}$ function in
norm. Fridly, Manchanda and Siddiqi \cite{fms} improve and extend the
results of Móricz and Siddiqi \cite{Mor} among them in $H_{p}$ norm, where $%
0<p<1.$ In \cite{gog8} Goginava investigated the behavior of Cesáro means of
Walsh-Fourier series in detail.

The main aim of this paper is to find necessary and sufficient conditions
for the convergence of Fejér means in terms of the modulus of continuity on
the Hardy spaces $H_{p},$ when $0<p\leq 1/2.$ In particular, the following
is true:

\begin{theorem}
Let $0<p\leq 1/2,$ $f\in H_{p},$ $M_{N}<n\leq $ $M_{N+1}$ and 
\begin{equation}
\omega \left( \frac{1}{M_{N}},f\right) _{H_{p}}=o\left( \frac{1}{%
M_{N}^{1/p-2}N^{2\left[ 1/2+p\right] }}\right) ,\text{ as \ }N\rightarrow
\infty .  \label{cond}
\end{equation}%
Then 
\begin{equation*}
\left\Vert \sigma _{n}f-f\right\Vert _{p}\rightarrow 0,\text{ when }%
n\rightarrow \infty .
\end{equation*}
\end{theorem}

\begin{theorem}
a) Let $0<p<1/2$ and $M_{N}<n\leq M_{N+1}.$ There exists a martingale $f\in
H_{p}(G_{m})$\ \ for which 
\begin{equation}
\omega \left( \frac{1}{M_{N}},f\right) _{H_{p}}=O\left( \frac{1}{%
M_{N}^{1/p-2}}\right) ,\text{ \ as \ }N\rightarrow \infty  \label{cond2}
\end{equation}%
and 
\begin{equation*}
\left\Vert \sigma _{n}f-f\right\Vert _{L_{p,\infty }}\nrightarrow 0,\,\,\,%
\text{as\thinspace \thinspace \thinspace }n\rightarrow \infty .
\end{equation*}%
b) Let $M_{N}<n\leq M_{N+1}.$ There exists a martingale $f\in H_{1/2}(G_{m})$%
\ \ for which 
\begin{equation}
\omega \left( \frac{1}{M_{N}},f\right) _{H_{1/2}}=O\left( \frac{1}{N^{2}}%
\right) ,\text{ \ as \ }N\rightarrow \infty  \label{cond3}
\end{equation}%
and 
\begin{equation*}
\left\Vert \sigma _{n}f-f\right\Vert _{1/2}\nrightarrow 0,\,\,\,\text{%
as\thinspace \thinspace \thinspace }n\rightarrow \infty .
\end{equation*}
\end{theorem}

\textbf{Proof of Theorem 1. }Let $f\in H_{p},$ $0<p\leq 1/2$ and $%
M_{N}<n\leq M_{N+1}.$ Using (\ref{1}) we have

\begin{eqnarray*}
&&\left\Vert \sigma _{n}f-f\right\Vert _{p}^{p} \\
&\leq &\left\Vert \sigma _{n}f-\sigma _{n}S_{M_{N}}f\right\Vert
_{p}^{p}+\left\Vert \sigma _{n}S_{M_{N}}f-S_{M_{N}}f\right\Vert
_{p}^{p}+\left\Vert S_{M_{N}}f-f\right\Vert _{p}^{p} \\
&=&\left\Vert \sigma _{n}\left( S_{M_{N}}f-f\right) \right\Vert
_{p}^{p}+\left\Vert S_{M_{N}}f-f\right\Vert _{p}^{p}+\left\Vert \sigma
_{n}S_{M_{N}}f-S_{M_{N}}f\right\Vert _{p}^{p} \\
&\leq &c_{p}\left( n^{1-2p}\log ^{2p\left[ 1/2+p\right] }n+1\right) \omega
^{p}\left( \frac{1}{M_{N}},f\right) _{H_{p}}+\left\Vert \sigma
_{n}S_{M_{N}}f-S_{M_{N}}f\right\Vert _{p}^{p}.
\end{eqnarray*}%
By simple calculation we get%
\begin{equation*}
\sigma _{n}S_{M_{N}}f-S_{M_{N}}f=\frac{M_{N}}{n}\left( \sigma
_{M_{N}}S_{M_{N}}f-S_{M_{N}}f\right) .
\end{equation*}

Let $p>0$ and $f\in H_{p}.$ Since (see \cite{We1} (for bounded Vilenkin
systems)) 
\begin{equation}
\Vert \sigma _{M_{k}}f-f\Vert _{p}\rightarrow 0,\text{ \ when }k\rightarrow
\infty  \label{31}
\end{equation}%
$S_{M_{N}}\sigma _{M_{N}}f=\sigma _{M_{N}}S_{M_{N}}f$ and $S_{M_{N}}f$ is
martingale, we obtain%
\begin{equation}
\left\Vert \sigma _{n}S_{M_{N}}f-S_{M_{N}}f\right\Vert _{p}^{p}=\frac{%
M_{N}^{p}}{n^{p}}\left\Vert \sigma _{M_{N}}S_{M_{N}}f-S_{M_{N}}f\right\Vert
_{p}^{p}\rightarrow 0,\text{ when }k\rightarrow \infty \text{.}  \label{10}
\end{equation}

It follows that under condition (\ref{cond}) we have%
\begin{equation*}
\left\Vert \sigma _{n}f-f\right\Vert _{p}\rightarrow 0,\text{ when }%
n\rightarrow \infty .
\end{equation*}%
Which complete the proof of Theorem 1.

\textbf{Proof of Theorem 2. } At first we consider case $0<p<1/2.$ Let 
\begin{equation*}
a_{k}=\frac{M_{k}^{1/p-1}}{\lambda }\left( D_{M_{k}+1}-D_{M_{k}}\right) ,%
\text{ }
\end{equation*}%
where $\lambda =\sup_{n\in \mathbb{N}}m_{n}$. It is easy to show that $a_{k}$
are $p$-atoms.

We set 
\begin{equation*}
f^{\left( A\right) }=\underset{i=0}{\overset{A}{\sum }}\frac{\lambda }{%
M_{i}^{1/p-2}}a_{i}.
\end{equation*}

Using Theorem W (see \cite{tep1}) we conclude that martingale $f\in H_{p}.$

On the other hand 
\begin{eqnarray}
&&f-S_{M_{n}}f  \label{20} \\
&=&\left( f^{\left( 1\right) }-S_{M_{n}}f^{\left( 1\right) },...,f^{\left(
n\right) }-S_{M_{n}}f^{\left( n\right) },...,f^{\left( n+k\right)
}-S_{M_{n}}f^{\left( n+k\right) }\right)  \notag \\
&=&\left( 0,...,0,f^{\left( n+1\right) }-f^{\left( n\right) },...,f^{\left(
n+k\right) }-f^{\left( n\right) },...\right)  \notag \\
&=&\left( 0,...,0,\underset{i=n+1}{\overset{n+k}{\sum }}\frac{a_{i}}{%
M_{i}^{1/p-2}},...\right) ,\text{ \ }k\in \mathbb{N}_{+}  \notag
\end{eqnarray}

is martingale. Applying Theorem W and (\ref{20}) we get%
\begin{equation*}
\omega (\frac{1}{M_{n}},f)_{H_{p}}\leq \sum\limits_{i=n+1}^{\infty }\frac{1}{%
M_{i}^{1/p-2}}\leq \frac{c}{M_{n}^{1/p-2}}=O\left( \frac{1}{M_{n}^{1/p-2}}%
\right) .
\end{equation*}

By simple calculation we have

\begin{equation}
\widehat{f}(j)=\left\{ 
\begin{array}{l}
M_{i},\,\ \ \text{if \thinspace \thinspace }j\in \left\{
M_{i},...,M_{i+1}-1\right\} ,\text{ }i=0,1,... \\ 
0,\text{ \thinspace \thinspace \thinspace\ \ if \thinspace \thinspace
\thinspace }j\notin \bigcup\limits_{i=0}^{\infty }\left\{
M_{i},...,M_{i+1}-1\right\} .\text{ }%
\end{array}%
\right.  \label{29}
\end{equation}

Combining (\ref{31}) and (\ref{29}) we can write%
\begin{eqnarray*}
&&\limsup\limits_{k\rightarrow \infty }\Vert \sigma _{M_{k}+1}f-f\Vert
_{L_{p,\infty }} \\
&=&\limsup\limits_{k\rightarrow \infty }\Vert \frac{M_{k}\sigma _{M_{k}}f}{%
M_{k}+1}+\frac{S_{M_{k}}(f)}{M_{k}+1}+\frac{M_{k}w_{M_{k}+1}}{M_{k}+1}-\frac{%
M_{k}f}{M_{k}+1}-\frac{f}{M_{k}+1}\Vert _{L_{p,\infty }} \\
&\geq &\limsup\limits_{k\rightarrow \infty }\frac{M_{k}}{M_{k}+1}\Vert
w_{M_{k}}\Vert _{L_{p_{\infty }}} \\
&&-\limsup\limits_{k\rightarrow \infty }\frac{M_{k}}{M_{k}+1}\Vert \sigma
_{M_{k}}f-f\Vert _{L_{p,\infty }}-\limsup\limits_{k\rightarrow \infty }\frac{%
1}{M_{k}+1}\Vert S_{M_{k}}f-f\Vert _{L_{p,\infty }} \\
&\geq &\limsup\limits_{k\rightarrow \infty }\frac{M_{k}}{M_{k}+1}\geq c>0.
\end{eqnarray*}

Now prove second part of Theorem 2. Let 
\begin{equation*}
f^{\left( A\right) }=\overset{A}{\sum_{i=1}}\frac{\lambda }{M_{i}^{2}}a_{i},
\end{equation*}%
where%
\begin{equation*}
a_{i}=\frac{M_{2M_{i}}}{\lambda }\left(
D_{M_{2M_{i}}+1}-D_{M_{2M_{i}}}\right)
\end{equation*}%
and $\lambda =\sup_{n\in \mathbb{N}}m_{n}$.

Analogously, we conclude that $f\in H_{1/2}$ and

\begin{equation*}
\omega (\frac{1}{M_{n}},f)_{H_{1/2}}\leq \overset{A}{\sum_{i=\left[ \log n/2%
\right] }}\frac{1}{M_{i}^{2}}=O\left( \frac{1}{n^{2}}\right) .
\end{equation*}%
Hence%
\begin{equation}
\sigma _{q_{M_{k}}}f-f=\frac{M_{2M_{k}}\sigma _{M_{2M_{k}}}(f)}{q_{M_{k}}}+%
\frac{1}{q_{M_{k}}}\sum_{j=M_{2M_{k}}+1}^{q_{M_{k}}}S_{j}f-\frac{M_{2M_{k}}f%
}{q_{M_{k}}}-\frac{q_{M_{k}-1}f}{q_{M_{k}}}  \label{nn}
\end{equation}%
It is easy to show that

\begin{equation}
\widehat{f}(j)=\left\{ 
\begin{array}{l}
\frac{M_{2M_{i}}}{M_{i}^{2}},\,\ \ \text{if \thinspace \thinspace }j\in
\left\{ M_{2M_{i}},...,M_{2M_{i}+1}-1\right\} ,\text{ }i=0,1,... \\ 
0,\text{ \ \ \ \ \ \thinspace \thinspace \thinspace if \thinspace \thinspace
\thinspace }j\notin \bigcup\limits_{i=0}^{\infty }\left\{
M_{2M_{i}},...,M_{2M_{i}+1}-1\right\} .\text{ }%
\end{array}%
\right.  \label{35}
\end{equation}%
Let $M_{2M_{k}}<j\leq q_{M_{k}}.$ Since $%
D_{j+M_{2M_{k}}}=D_{M_{2M_{k}}}+w_{M_{2M_{k}}}D_{j},$ \ when \thinspace
\thinspace $j<M_{2M_{k}}$ \ using (\ref{35}) we have 
\begin{equation*}
S_{j}f=S_{M_{2M_{k}}}f+\sum_{v=M_{_{2M_{k}}}}^{j-1}\widehat{f}%
(v)w_{v}=S_{M_{2M_{k}}}f+\frac{%
M_{2M_{k}}w_{M_{_{2M_{k}}}}D_{_{j-M_{_{2M_{k}}}}}}{M_{k}^{2}}.
\end{equation*}

and 
\begin{eqnarray*}
&&\frac{1}{q_{M_{k}}}\sum_{j=M_{_{2M_{k}}}+1}^{q_{M_{k}}}S_{j}f=\frac{%
q_{M_{k}-1}S_{M_{M_{k}}}f}{q_{M_{k}}}+\frac{M_{2M_{k}}w_{M_{_{2M_{k}}}}}{%
q_{M_{k}}M_{k}^{2}}\sum_{j=1}^{q_{M_{k}-1}}D_{_{j}} \\
&=&\frac{q_{M_{k}-1}S_{M_{2M_{k}}}f}{q_{M_{k}}}+\frac{%
M_{2M_{k}}w_{M_{_{2M_{k}}}}q_{M_{k}-1}K_{q_{M_{k}-1}}}{q_{M_{k}}M_{k}^{2}}.
\end{eqnarray*}

Using (\ref{nn}) we get%
\begin{eqnarray}
&&\Vert \sigma _{q_{M_{k}}}f-f\Vert _{1/2}^{1/2}\geq \frac{c}{M_{k}}\Vert
q_{M_{k}-1}K_{q_{M_{k}-1}}\Vert _{1/2}^{1/2}  \label{a11} \\
&&-\left( \frac{M_{2M_{k}}}{q_{M_{k}}}\right) ^{1/2}\Vert \sigma
_{M_{_{2M_{k}}}}f-f\Vert _{1/2}^{1/2}-\left( \frac{q_{M_{k}-1}}{q_{M_{k}}}%
\right) ^{1/2}\Vert S_{M_{2M_{k}}}f-f\Vert _{1/2}^{1/2}.  \notag
\end{eqnarray}

Let $x\in I_{2s+1}\left( s_{2k}e_{2k}+s_{2\eta }e_{2\eta }\right) ,$ $1\leq
s_{2\eta }\leq m_{2\eta }-1$ and $1\leq s_{2k}\leq m_{2k}-1,$ $\eta
=1,...,M_{k}-4,\,$\ $s=\eta +2,\eta +3,...,M_{k}-2.$ Applying (\ref{3a}) we
have

\begin{equation*}
q_{M_{k}-1}\left\vert K_{q_{M_{k}-1}}\left( x\right) \right\vert \geq \frac{%
M_{2\eta }M_{2s}}{4}.
\end{equation*}

We can write

\begin{eqnarray}
&&\int_{G_{m}}\left\vert q_{M_{k}-1}K_{q_{M_{k}-1}}\right\vert ^{1/2}d\mu
\label{33} \\
&\geq &\frac{1}{4}\sum_{\eta =1}^{M_{k}-4}\sum_{s=\eta
+2}^{M_{k}-2}\int_{I_{2s+1}\left( s_{2k}e_{2k}+s_{2l}e_{2l}\right)
}\left\vert q_{M_{k}-1}K_{q_{M_{k}-1}}\right\vert ^{1/2}d\mu  \notag \\
&\geq &\frac{1}{4}\sum_{\eta =1}^{M_{k}-4}\sum_{s=\eta +2}^{M_{k}-2}\frac{1}{%
2M_{2s}}\sqrt{M_{2s}M_{2\eta }}\geq cM_{k}.  \notag
\end{eqnarray}

Combining (\ref{31}), (\ref{a11}-\ref{33}) we have%
\begin{equation*}
\underset{k\rightarrow \infty }{\lim \sup }\Vert \sigma _{q_{M_{k}}}f-f\Vert
_{1/2}\geq c>0.
\end{equation*}

Theorem 2 is proved.

\end{document}